\documentclass[11pt]{amsproc}
\usepackage{amsmath}
\usepackage{amssymb}
\usepackage{amsfonts}

\newtheorem{theorem}{Theorem}
\newtheorem{proposition}{Proposition}
\newtheorem{corollary}{Corollary}
\newtheorem{example}{Example}

\newtheorem{lemma}{Lemma}

\begin{document}
\title{Factorial threefolds with $\mathbb{G}_a$-actions}
\author{Adrien Dubouloz, David R. Finston, and Parag Deepak Mehta}
\address{Adrien Dubouloz\\
CNRS\\
Institut de Math\'{e}matiques de Bourgogne\\
Universit\'{e} de Bourgogne\\
9 Avenue Alain Savary\\
BP 47870\\
21078 Dijon Cedex\\
France}
\email{Adrien.Dubouloz@u-bourgogne.fr}
\address{David R. Finston\\
Department of Mathematical Sciences\\
New Mexico State University\\
Las Cruces, New Mexico 88003}
\email{dfinston@nmsu.edu}
\address{Parag Deepak Mehta\\
Mahindra United World College of India\\
Pune 412108 \\
India }
\email{pmehta@muwci.net }

\begin{abstract}
The affine cancellation problem, which asks whether complex affine varieties
with isomorphic cylinders are themselves isomorphic, has a positive solution
for two dimensional varieties whose coordinate rings are unique
factorization domains, in particular for $\mathbb{C}^{2}$,\ but
counterexamples are found within normal surfaces (Danielewski surfaces) and
factorial threefolds of logarithmic Kodaira dimension equal to $1$. \ \ The
latter are therefore remote from $\mathbb{C}^{3}$, the first unknown case
where the base of one cylinder is an affine space. Locally trivial $\mathbb{G%
}_{a}$-actions play a significant role in these examples.\ Factorial
threefolds admitting free $\mathbb{G}_{a}$-actions are discussed, especially
a class of varieties with negative logarithmic Kodaira dimension which are
total spaces of nonisomorphic $\mathbb{G}_{a}$-bundles. \ Some members of
the class are shown to be isomorphic as abstract varieties,\ but it is
unknown whether any members of the class constitute counterexamples to
cancellation. \ 
\end{abstract}

\maketitle

\section{Introduction}

Let $\mathbb{G}_a$ denote the additive group of complex numbers, and $X$ a
complex affine variety. Throughout, an action of $\mathbb{G}_a$ on $X$ means
an algebraic action. Every such action can be realized as the exponential of
some locally nilpotent derivation $\delta $ of the coordinate ring $\mathbb{C%
}[X]$ of $X$ and every locally nilpotent derivation gives rise to an action.
The ring $C_{0}$ of $\mathbb{G}_a$-invariants in $\mathbb{C}[X]$ is equal to
the ring of constants of the generating derivation.

An action $\sigma :\mathbb{G}_a\times X\rightarrow X$ is said to be
equivariantly trivial if there is a variety $Y$ for which $X$ is $\mathbb{G}%
_a$-equivariantly isomorphic to $Y\times \mathbb{G}_a,$ where $\mathbb{G}_a$
acts on $Y\times \mathbb{G}_a$ by translations on the second factor. In this
case $Y$ is affine and $X$ the total space of a trivial $\mathbb{G}_a$%
-bundle over $Y.$ A section for the bundle can be identified with the zero
locus of a regular function $s\in \mathbb{C}[X]$ for which $\delta s=1.$
Such a function is called a slice; if one exists, $\mathbb{C}[X]=C_{0}[s]$
and $Y\cong $ \textrm{\ Spec }$C_{0}.$ The action is said to be locally
trivial if a geometric quotient $\pi :X\rightarrow Y$ exists in the category
of algebraic spaces for which $\pi $ is a principal $\mathbb{G}_a$-bundle.
The action is said to be proper if the morphism 
\begin{equation*}
\sigma \times id:\mathbb{G}_a\times X\rightarrow X\times X
\end{equation*}
is proper.

The main point of the paper is to investigate the class of hypersurfaces 
\begin{equation*}
X_{m,n}:x^{m}v-y^{n}u=1
\end{equation*}%
in $\mathbb{A}^{4}=\mathrm{Spec}(\mathbb{C}[x,y,u,v])$, where $m,\, n$ are
positive integers. \ Since the polynomial $x^{m}v-y^{n}u=1$ is an invariant
for the $\mathbb{G}_a$ action 
\begin{equation*}
t(x,y,u,v)=(x,y,u+tx^{m},v+ty^{n})
\end{equation*}
on $\mathbb{A}^{4},$ the $X_{m,n}$ are stable under this action. \ Moreover,
with $\mathbb{A}^2_*= \mathbb{A}^{2}\setminus\{0\},$ the complement of the
origin in the complex affine plane, the $X_{m,n}$, or rather the cohomology
classes in $H^{1}(\mathbb{A}^2_*,O_{\mathbb{A}^2_*})$ of the \v{C}ech
coycles $x^{-m}y^{-n}\in\mathbb{C}[x^{\pm 1}, y^{\pm 1}]$, constitute a
natural basis for the set of isomorphy classes of principal $\mathbb{G}_a$%
-bundles over $\mathbb{A}^2_*$. As seen below, these varieties can be viewed
as globalizing the local structure of arbitrary smooth threefolds with
locally trivial $\mathbb{G}_a$-action.

Similarly constructed threefolds were investigated in \cite{F-M} as
factorial counterexamples to a generalized cancellation problem: 
\begin{equation*}
\ \text{If }X\times \mathbb{A}^1\cong Y\times \mathbb{A}^1\text{ is }X\
\cong Y \text{ ?}
\end{equation*}
Those examples arise from principal $\mathbb{G}_a$-bundles over the smooth
locus of factorial surfaces of logarithmic Kodaira dimension 1, hence admit
no nonconstant morphism from $\mathbb{A}^{2}$. As such they differ
substantially from affine three space. As total spaces for $\mathbb{G}_a$%
-bundles over $\mathbb{A}^{2}_*$ the varieties $X_{m,n}$ considered here are
closer in spirit to affine spaces and also have isomorphic cylinders. \ For
distinct pairs $(m,n)$ the corresponding varieties $X_{m,n}$ represent
nonisomorphic $\mathbb{G}_a$-bundles, nevertheless some of them turn out to
be isomorphic as varieties. \ \ We do not know as yet whether any\ pair of
them provides a counterexample to cancellation.

For factorial affine varieties, i.e. those whose coordinate ring is a unique
factorization domain (UFD), local triviality is equivalent with the
intersection of $C_{0}$ with the image of $\delta $ generating the unit
ideal in $\mathbb{C}[X]$ \cite{D-F}$.$ \ \ In this case the geometric
quotient has the structure of a quasiaffine variety. To see this explicitly,
consider a locally trivial $\mathbb{G}_a$-action on a factorial affine
variety $X$ and let $\delta (a_{1}),...,\delta (a_{n})\in C_{0}$ generate
the unit ideal in $\mathbb{C}[X].$ The geometric quotient is embedded as an
open subvariety in \textrm{Spec }$R$ for an integrally closed ring $R$
constructed as follows. Set $R_{i}=\mathbb{C}[X,\frac{1}{\delta (a_{i})}]^{%
\mathbb{G}_a}.$ \ Note that $\mathbb{C}[X,\frac{1}{\delta (a_{i})}]=R_{i}[%
\frac{a_{i}}{\delta (a_{i})}]$ so that $R_{i}$ is a finitely generated $%
\mathbb{C}$-algebra, say 
\begin{equation*}
R_{i}=\mathbb{C}[b_{i1},...,b_{im},\frac{1}{\delta (a_{i})}],
\end{equation*}%
with $b_{ij}\in C_{0.}$ \ Define $R$ to be the integral closure of $\mathbb{C%
}\mathbf{[}b_{ij},\delta (a_{i})|1\leq i\leq n,1\leq j\leq m].$ \ The map $%
\pi :X\rightarrow \mathrm{Spec}(R)$ induced from the ring inclusion is
smooth, hence open,\ and endows its image\ $U$ with the structure of a
geometric quotient. \ Moreover, for the open cover $\mathcal{U}=\{U_{i}=%
\mathrm{Spec}(R_{i})\}_{1\leq i \leq n}$ of $U$, the \v{C}ech cocycle 
\begin{equation*}
(\frac{a_{1}}{\delta (a_{1})}-\frac{a_{2}}{\delta (a_{2})},\ldots ,\frac{%
a_{n-1}}{\delta (a_{n-1})}-\frac{a_{n}}{\delta (a_{n})})\in C^{1}(\mathcal{U}%
,O_{U})
\end{equation*}%
defines the structure of $\pi :X\rightarrow U$ as a principal $\mathbb{G}_a$%
-bundle.

Recall that for an integral domain $R$ with quotient field $K,$ the
transform of $R$ with respect to the ideal $I$, written $T_{I}(R)$ is $\cup
_{n>0}\{\alpha \in K|\alpha I^{n}\subset R\}$. The ring of $\mathbb{G}_a$%
-invariants is isomorphic to the ring of sections of the structure sheaf of $%
\mathrm{Spec }(R)$ on $U$, hence to $T_{I}(R)$ where $I$ defines the closed
subscheme $Z=\mathrm{Spec}(R)\setminus U$ equipped with its reduced
structure. It is clear from the construction that $I=\sqrt{(\delta
(a_{1}),...,\delta (a_{n}))}$. Since $R$ is normal, $R=C_{0}$ if all
components of $Z$ have codimension $\geq 2$ (equivalently if $I$ is
contained in no height one prime ideal of $R$).

In case the dimension of $X$ is less than or equal to three, a theorem of
Zariski \cite[p.45]{Nag} implies that $C_{0}$ is a finitely generated $%
\mathbb{C}$-algebra. Combined with the fact that the complement of a pure
codimension one subvariety of a two dimensional normal affine variety is
again an affine variety \cite[p.45]{Nag}, the following ring theoretic
criterion for triviality of a $\mathbb{G}_a$-bundle with a three dimensional
affine factorial total space $X$ is obtained: \indent\newline

\emph{Let }$X$\emph{\ be an affine factorial threefold with a locally
trivial }$\mathbb{G}_a$\emph{-action generated by the locally nilpotent
derivation }$\delta $\emph{\ of }$C[X]$\emph{. \ With the ring }$R$\emph{,
ideal }$I,$\emph{and }$Z,U$\emph{\ as above, let }$J=\cap \underset{ht\text{ 
}\mathfrak{p=1}}{_{I\subset \mathfrak{p}}}p$\emph{\ . Then }$C_{0}$\emph{\ }$%
=T_{J}R$\emph{\ and the action is equivariantly trivial if and only if }$J=I$%
\emph{. } \newline

Indeed, since the complement\ $W$ of the zero locus of $J$ in $\mathrm{Spec}%
( R)$ is affine, every principal $\mathbb{G}_a$-bundle over $W$ is trivial.
\ If $J$ properly contains $I$ then $U$ is equal to the complement of a
finite but nonempty subset of $\mathrm{Spec}(C_{0})$ (since $I$ is\ a
radical ideal, the Nullstellensatz implies that it has no embedded primes).
\ In particular, $U$ is not affine, and therefore neither is $U\times 
\mathbb{A}^1.$

More generally, a set-theoretically free action of an algebraic group on a
normal variety admits a geometric quotient in the category of algebraic
spaces, and the quotient map is a locally trivial fiber bundle in the \'{e}%
tale topology. In this context, the properness of the action is then
equivalent to the separatedness of the algebraic space quotient. Restricting
now to proper $\mathbb{G}_{a}$-actions on normal quasiaffine threefolds, we
obtain that the quotient exists as a quasiprojective variety and that the
quotient map is locally trivial in the Zariski topology. Indeed, the
quotient embeds as an open subset of a normal two dimensional analytic
space. By Chow's lemma we know that a desingularization is quasiprojective,
hence the quotient is the image of a quasiprojective variety under a proper
morphism. \ As such it is again quasiprojective. \ That the quotient
morphism is a Zariski locally trivial $\mathbb{G}_{a}$-bundle follows for
instance from \cite[Prop. 0.9]{Mum}

All proper $\mathbb{G}_a$-actions on quasiaffine surfaces are known to be
equivariantly trivial. \ For normal surfaces this was observed in \cite{D-F2}%
, but again, an application of \cite[Prop. 0.9]{Mum} gives the result in
general. All fixed point free actions on contractible threefolds are
equivariantly trivial \cite{K-S}. Examples of a nonlocally trivial proper
action on $\mathbb{A}^{5}$ and on a smooth factorial fourfold are also
discussed in \cite{D-F2}. The latter example is isomorphic to a cylinder
over an affine variety which itself has the structure of the total space of
a principal $\mathbb{G}_a$-bundle over the punctured affine plane $\mathbb{A}%
^{2}_*$, the topic of the subsequent sections of the present work. Except
for special kinds of actions \cite{D-F3}, the issue of local triviality of
proper $\mathbb{G}_a$-actions on the affine four-space $\mathbb{A}^{4}$ is
unsettled.

Replacing the smoothness hypothesis on the affine threefold $X$ by
factoriality and smoothness of the ring of $\mathbb{G}_a$-invariants, the
following lemma of Miyanishi \cite{miya} leads to a local description of $X$.

\begin{lemma}
Let $(\mathcal{O},\mathfrak{m})$ be a regular local ring of dimension $n\geq
2$ and let $A$ be a factorial, finitely generated $\mathcal{O}$-domain with $%
\mathcal{O}\hookrightarrow A.$ Let $f:X\rightarrow Y$ be the morphism
induced by the ring inclusion, where $X= \mathrm{Spec}(A)$ and $Y=\mathrm{%
Spec}(\mathcal{O}).$ Let $U=Y\setminus\!\{ \mathfrak{m}\}.$ Assume that $%
f_{U}:f^{-1}(U)\rightarrow U$ is an $\mathbb{A}^{1}$-bundle. Then either $%
X\cong Y\times \mathbb{A}^{1}$ or $f^{-1}(\mathfrak{m})=\emptyset $ (the
latter is only possible if $n=2$).\ 
\end{lemma}

\begin{corollary}
Let $X$ be a smooth affine factorial threefold with a locally trivial $%
\mathbb{G}_a$-action that is not equivariantly trivial. \ Let $\delta $ be
the derivation generating the action and $I=\sqrt{C_{0}\cap \ker \delta }%
\mathbb{C}[X]$. \ Assume that $C_{0}$ is regular. \ For a maximal ideal $%
\mathfrak{m}$ of $C_{0},$ set $S=C_{0}\setminus\mathfrak{m}$. \ Then with $%
u,v$ denoting algebraically independent indeterminates,

\begin{enumerate}
\item $S^{-1}\mathbb{C[}X\mathbb{]\cong }(C_0)_{\mathfrak{m}}[u]$ if $%
I\nsubseteqq \mathfrak{m}$ and

\item $S^{-1}\mathbb{C[}X\mathbb{]\cong }(C_0)_{\mathfrak{m}}[u,v]/(au-bv-1)$
for some regular sequence $(a,b)$ in $\mathfrak{m}$ if $I\subseteq \mathfrak{%
m.}$
\end{enumerate}

\begin{proof}
Since $S\subset C_{0},$ the action extends to $S^{-1}\mathbb{C[}X\mathbb{]}.$
\ Since the action is not equivariantly trivial , $IC_{0}\neq C_{0}.$ If $%
I\nsubseteqq \mathfrak{m,}$ then $S^{-1}\mathbb{C[}X\mathbb{]}$ contains a
slice. \ If $I\subseteq \mathfrak{m,}$ then apply the lemma to $\mathcal{O}%
=(C_0)_{\mathfrak{m}}$ and $A=$ $S^{-1}\mathbb{C[}X]$ to see that $\mathfrak{%
m}S^{-1} \mathbb{C}[X]=C[X].$ On the other hand, no proper principal ideal
of $(C_0)_{\mathfrak{m}}$ blows up in $S^{-1}\mathbb{C}[X].$
\end{proof}
\end{corollary}

The $X_{m,n}$ above represent a class of global versions of case 2.

\begin{example}
Using a cover by the principal open subsets $x\neq 0,y\neq 0,$ one sees that
a \textquotedblleft basis" for the principal $\mathbb{G}_{a}$-bundles over $%
\mathbb{A}_{\ast }^{2}$ is given by the affine varieties $X_{m,n}\subset 
\mathbb{A}^{4}$ defined by $x^{m}v-y^{n}u-1=0$ where $m,\,n$ are positive
integers. \ These correspond to the \v{C}ech 1-cocycles $x^{-m}y^{-n}$
relative to the given cover of $\mathbb{A}_{\ast }^{2}$ . The $X_{m,n}$ are
smooth factorial threefolds with locally trivial $\mathbb{G}_{a}$-action
generated by the locally nilpotent derivation $\delta :u\mapsto x^{m}\mapsto
0,$ \ $v\mapsto y^{n}\mapsto 0.$ \ Since $\mathbb{A}_{\ast }^{2}$ is a
geometric quotient, $C_{0}=\mathbb{C}[x,y]\cong \Gamma (\mathbb{A}_{\ast
}^{2},\mathcal{O}_{\mathbb{A}_{\ast }^{2}}).$
\end{example}

Regularity of $C_{0}$\ is not a necessary condition for the conclusion of
the corollary. \ 

\begin{example}
Let $A=\mathbb{C[}x,y,z]/(x^{2}+y^{3}+z^{5}),$ the coordinate ring for the
analytically factorial singular surface $S,$ and set $U=S\setminus \{0\},$
the smooth locus of $S$ $.$ Using the cover of $U$ by the two open subsets $%
S_{x}\cup S_{y},$, the cohomology classes of the \v{C}ech cocycles $%
x^{-m}y^{-n}z^{k}$, $0\leq k\leq 4$ constitute a basis for $H^{1}(U,O_{U})$.
Consider the affine threefolds $Z_{m,n,k}\subset \mathbb{A}^{5}$ defined by 
\begin{equation*}
x^{m}v-y^{n}u-z^{k}=0,\;x^{2}+y^{3}+z^{5}=0,0\leq k\leq 4
\end{equation*}%
for $0\leq k\leq 4$ $.$ \ With $\pi _{m,n,k}:Z_{m,n,k}\rightarrow S$ induced
by the ring inclusion $A\subset \mathbb{C[}Z_{m,n,k}\mathbb{]}$, \ the
cocycle $x^{-m}y^{-n}z^{k}$ corresponds to the principal $\mathbb{G}_{a}$%
-bundle over $U$ with total space $E_{m,n,k}=\pi ^{-1}(U)$ and bundle
projection $\pi $. Then $E_{m,n,k}$ is smooth, factorial, affine threefold
admitting a locally trivial $\mathbb{G}_{a}$-action with quotient isomorphic
to $U$.
\end{example}

Smoothness of $E_{m,n,k}$ is clear from the defnition and factorialty
follows from the fact that $\pi$ is a Zariski fibration with base and fiber
having trivial Picard group \cite{Magid}. \ In fact $E_{m,n,k}$ is the
quotient of a principal $\mathbb{G}_a$-bundle over $\mathbb{A}^{2}_*$ by the
action of a finite group (the binary icosahedral group). \ We will see in
the next section that all nontrivial $\mathbb{G}_a$-bundles over $\mathbb{A}%
^{2}_*$ have affine total spaces and deduce that the $E_{m,n,k}$ are affine
by Chevalley's theorem.

The counterexamples to cancellation in \cite{F-M} are constructed
analogously to the $E_{m,n,k},$ i.e$.$ as principal $\mathbb{G}_a$-bundles
over surfaces 
\begin{equation*}
S_{a,b,c}:x^{a}+y^{b}+z^{c}=0,\text{ }\frac{1}{a}+\frac{1}{b}+\frac{1}{c}<1
\end{equation*}%
punctured at the singular point. The numerical criterion $\frac{1}{a}+\frac{1%
}{b}+\frac{1}{c}<1$ enabled the distinction of the total spaces as affine
varieties. \ We do not have analogous results for the $E_{m,n,k}.$

\section{Principal $\mathbb{G}_a$-bundles over $\mathbb{A}^{2}_*=\mathbb{A}%
^{2}\setminus\{0\}$}

We start with the observation that the total spaces of nontrivial principal $%
\mathbb{G}_{a}$-bundles over $\mathbb{A}_{\ast }^{2}$ are always affine
schemes. Indeed, a nontrivial \v{C}ech cocycle with value in $\mathcal{O}_{%
\mathbb{A}_{\ast }^{2}}$ for the covering of $\mathbb{A}_{\ast }^{2}$ by the
open subsets $\{x\neq 0\}$ and $\{y\neq 0\}$ can be written as $%
g=p(x,y)x^{-m}y^{-n}\in\mathbb{C}[x^{\pm 1},y^{\pm 1}]$ with $\deg
_{x}p<m,\, \deg _{y}p<n$ and denote by $X(m,n,p)$ the total space for the so
determined bundle over $\mathbb{A}_{\ast }^{2}$. Set 
\begin{equation*}
A=\mathbb{C}[x,y,u,v]/(x^{n}v-y^{m}u-p(x,y))
\end{equation*}%
so that $X(m,n,p)$ is isomorphic to $\mathrm{Spec}(A)\setminus \!V(I)$ where 
$I=(x,y)A$.

If $p(0,0)\neq 0$ then $X(m,n,p)$ is the zero locus of $x^{m}v-y^{n}u-p(x,y)$
in $\mathbb{A}^{4}.$ \ Assume then that $p(0,0)=0$ and note that $I$ is a
height one prime ideal in $A.$ \ We show that $IT_{I}(A)=T_{I}(A)\ $to
conclude that $X(m,n,p)=\mathrm{Spec}(A)\setminus V(I)$ is affine. Write $%
p(x,y)=\sum_{i=0}^{k}p_{i}(x)y^{i}$

\emph{Case 1}: $p_{0}(x)\neq 0.$ \ Let $a\geq 1$ be the multiplicity of $0$
as a root of $p$ and write $p_{0}(x)=x^{a}q_{0}(x).$ From 
\begin{equation*}
p(x,y)=x^{a}q_{0}(x)+y\sum_{i=1}^{k}p_{i}(x)y^{i-1}
\end{equation*}%
we obtain 
\begin{equation*}
\frac{x^{m-a}-q_{0}}{y}=\frac{y^{n-1}u+\sum_{i=1}^{k}p_{i}(x)y^{i-1}}{x^{a}}%
\in T_{I}(A).
\end{equation*}%
Thus $q_{0}(x)\in IT_{I}(A),$ but $q_{0}(0)\neq 0$ implies that $%
IT_{I}(A)=T_{I}(A).$

\emph{Case} \emph{2}: $\ p_{0}(x)=0.$ Write $p(x,y)=y^{b}%
\sum_{i=0}^{k-b}q_{i}(x)y^{i}$ with $q_{0}(x)\neq 0.$ From 
\begin{equation*}
x^{m}v=y^{b}[y^{n-b}u+\sum_{i=0}^{k-b}q_{i}(x)y^{i}]
\end{equation*}%
we obtain%
\begin{equation*}
w=\frac{v}{y^{b}}=\frac{y^{n-b}u+\sum_{i=0}^{k-b}q_{i}(x)y^{i}]}{x^{m}}\in
T_{I}(A).
\end{equation*}%
Then $v=y^{b}w$ so that in $T_{I}(A)$,%
\begin{eqnarray*}
x^{m}y^{b}w-y^{n}u &=&y^{b}\sum_{i=0}^{k-b}q_{i}(x)y^{i} \\
x^{m}w-y^{n-b}u &=&\sum_{i=0}^{k-b}q_{i}(x)y^{i}.
\end{eqnarray*}

Replace $A$ with $A^{\prime }=\mathbb{C}[x,y,u,v]/(x^{m}w-y^{n-b}u-%
\sum_{i=0}^{k-b}q_{i}(x)y^{i})$ to reduce to Case 1. \newline

The topological universal cover of the $\ E_{m,n,k}$ of the previous section
can be viewed as a principal $\mathbb{G}_a$-bundle over $\mathbb{A}^{2}_*$
and therefore as a "linear combination" of $X_{m,n}$ in the sense of \v{C}%
ech 1-cocycles. \ Indeed, for $G$ the binary icosahedral group embedded in $%
GL_{2}(\mathbb{C)}$ it is well known that $\mathbb{A}^{2}_*/G\cong U$ and
that the \'{e}tale mapping $\pi :\mathbb{A}^{2}_*\rightarrow \mathbb{A}%
^{2}_*/G$ is the universal topological covering. \ Thus, given a nontrivial $%
\mathbb{G}_a $-bundle $E_{m,n,k}$ over $U,$ we obtain by base extension $%
\widetilde{E}_{m,n,k}=\mathbb{A}^{2}_*\times _{U}\ E_{m,n,k}$ a nontrivial
principal $\mathbb{G}_a$-bundle over $\mathbb{A}^{2}_*$ whose total space
yields the universal covering space of $\ E_{m,n,k}$. Moreover, $\ E_{m,n,k}$
is recovered as $\widetilde{E}_{m,n,k}/G$ where the finite group $G$ acts
freely on the first factor. It follows that each $\widetilde{E}_{m,n,k}$ is
affine, smooth, and factorial. \ Finally, since $\widetilde{E}%
_{m,n,k}\rightarrow \ E_{m,n,k}\cong \ \ \widetilde{E}_{m,n,k}/G$ is a
finite morphism, $\ E_{m,n,k}$ is affine (and smooth, factorial) by
Chevalley's theorem.

The varieties $X_{m,n}$ (resp. $E_{m,n,k}$) in the preceding examples are
total spaces for principal $\mathbb{G}_a$-bundles over the quasiaffine $%
\mathbb{A}^{2}_*$ (resp. $U).$ The base extension $X_{m,n}\times_{\mathbb{A}%
^{2}_*} X_{p,q}$ is therefore a principal $\mathbb{G}_a$-bundle over both $%
X_{m,n}$ and $X_{p,q}$ for any $m,\,n,\,p,\,q$. Since they are affine, we
obtain $X_{m,n}\times \mathbb{A}^1$ $\cong $ $X_{p,q}\times \mathbb{A}^1$
and the same holds true for the cylinders over the $\ E_{m,n,k}$.

The well known Danielewski surfaces $D(n)$ given as the zero loci of the
polynomials $y^{2}-2x^{n}z-1,$ $n>0$ in $\mathbb{A}^{3}$ provide smooth but
nonfactorial counterexamples to the generalized cancellation problem. \ In
fact for $n\neq m,$ $D(n)$ and $D(m)$ are not even homeomorphic in the
Euclidean topology. \ In contrast, all $X_{m,n}$ (resp. $\ E_{m,n,k})$ are
homeomorphic to $\mathbb{A}^2_*\times \mathbb{A}^1$ (resp. $U\times \mathbb{A%
}^1$ since, as principal bundles, the Euclidean topology of the base has a
countable basis and the fiber is solid \cite{Steenrod}.

These considerations apply to any of the affine surfaces in $\mathbb{A}^{3}$
with a rational double point (quotient singularity) at the origin. \
Principal $\mathbb{G}_a$-bundles over the quasiaffine surface obtained by
puncturing at the singular point will again produce potential
counterexamples to the generalized affine cancellation problem. The surface
defined by $x^{2}+y^{3}+z^{5}=0$ was chosen for its historical significance
and the fact that its singularity is the unique one that is analytically
factorial (hence "closest" to the affine plane $\mathbb{A}^{2}$). \ But we
do not know whether or not the total spaces are \ isomorphic as varieties

On the other hand, the construction does provide cancellation
counterexamples for other quotient singularities. \ The following example
relies on the Makar Limanov invariant of an affine domain $A$ : 
\begin{equation*}
ML(A)=\cap \{\ker \delta :\delta \text{ is a locally nilpotent derivation of 
}A\}
\end{equation*}

\begin{example}
Let $U=S\setminus\{0\}$ where $S\subset \mathbb{A}^{3}$ is defined by 
\begin{equation*}
x^{2}+y^{2}z+z^{c}=0,\text{ }c\geq 3.
\end{equation*}%
It is well known that $S\cong \mathbb{A}^{2}/G$ where $G$ is the binary
dihedral group of order $4c.$ Define $X_{m,n}(G)\subset S\times \mathbb{A}%
^{2}$ by $x^{m}v-y^{n}u-1=0$ where $m,\, n$ are positive integers. Then
\end{example}

\begin{enumerate}
\item $X_{m,n}(G)$ is not factorial, but $ML(\mathbb{C}[X_{m,n}(G)])=\mathbb{%
C}[S]$. To see this, let $\delta $ be locally nilpotent with kernel $R$, and 
$K=\mathrm{Frac}(R)$. Then the extension of $\delta $ to $K\otimes _{R}%
\mathbb{C}[X_{m,n}(G)]$ has a slice $s$ and 
\begin{equation*}
K\otimes _{R}\mathbb{C}[X_{m,n}(G)]\ \cong K[s],\text{ }\delta (s)=1
\end{equation*}%
and $K$ is the kernel of the extended derivation. \ If $\deg _{s}(z)>0,$
then $z(s)$ divides both $x(s),y(s)$ contradicting $x^{m}v-y^{n}u=1.$ Thus $%
z\in R$ from which it follows that both $x,y\in R.$

\item $X_{m,n}(G)\cong X_{p,q}(G)$ if and only if $(m,n)=(p,q)$. This
follows from arguments analogous to those in \cite{F-M} and the fact that
any automorphism of $\mathbb{C}[Y]$ is of the form $(x,y,z)\mapsto (\mu
_{1}x,\mu _{2}y,\mu _{3}z)$ for certain roots of unity $\mu _{i}$ \cite{M-M}.

\item $\mathbb{A}^{2}_*\times _{U}X_{m,n}(G)\rightarrow X_{m,n}(G)$ is an 
\'{e}tale covering.
\end{enumerate}

Note that the total spaces $\mathbb{A}^{2}_*\times _{U}X_{m,n}(G)$ of these
principal $\mathbb{G}_a$-bundles over $\mathbb{A}^{2}_*$ are not isomorphic
as $G$-varieties, but it is possible that there is some abstract isomorphism
between them.

It is a classical result going back to Klein that the ring of invariants $%
\mathbb{C}[x,y]^{G}$ for the natural action of a finite subgroup $G$ of $%
SL_{2}(\mathbb{C)}$ on $\mathbb{A}^{2}$ corresponding to the quotient
surface singularities can be generated by 3 homogeneous polynomials. As a
consequence, affine threefolds of the type $X_{m,n}(G)$ for such $G$ can be
interpreted as quasihomogeneous varieties with respect to a grading
determined by the generators of $\mathbb{C}[x,y]^{G}.$ Properties of such
quasihomogeneous threefolds will be addressed in future work.

\section{Isomorphy Classes of $X_{m,n}$}

Recall the hypersurfaces $X_{m,n}\subset \mathbb{A}^{4}=\mathrm{Spec}(%
\mathbb{C}[x,y,u,v])$ defined by 
\begin{equation*}
x^{m}v-y^{n}u=1
\end{equation*}
with $mn\neq 0.$ The main result of this section shows that at least for
some distinct pairs $\{m,n\}\neq \{p,q\}$ we have $X_{m,n}\cong X_{p,q}$.

In his investigation of complex surfaces with $\mathbb{G}_{a}$-actions \cite%
{tD}, tom Dieck observed that $X_{m,m}$ admits a free action of a semidirect
product structure $G_{2m}\ $ on $\mathbb{G}_{m}\times \mathbb{G}_{a}$ for
which $X_{m,m}/G_{2m}\cong \mathbb{P}^{1}$ and $X_{m,m}\rightarrow \mathbb{P}%
^{1}$ is a principal $G_{2m}$-bundle. Moreover, the induced mapping $X_{m,m}/%
\mathbb{G}_{m}\rightarrow \mathbb{P}^{1}$ is diffeomorphic to the complex
line bundle $\mathcal{O}_{\mathbb{P}^{1}}(-2m).$ The quotient $X_{m,m}/%
\mathbb{G}_{m}$ turns out to be isomorphic to the Danielewski surface $%
D(m):y^{2}-2x^{m}z=1$ and the latter can thus be distinguished topologically
for distinct values of $m$. A related construction enables us to prove the
existence of an algebraic isomorphism $X_{m,n}\cong X_{p,q}$ if $m+n=p+q$
albeit with a different analysis of the structure.

For a given $(m,n)$ \ we have actions of $\mathbb{G}_m$ and $\mathbb{G}_a$
on $X_{m,n}$ given respectively by 
\begin{eqnarray*}
\mathbb{G}_m\times X_{m,n} &\rightarrow &X_{m,n} \\
(\lambda ,(x,y,u,v)) &\longmapsto &(\lambda x,\lambda y,\lambda
^{-n}u,\lambda ^{-m}v).
\end{eqnarray*}
and 
\begin{eqnarray*}
\mathbb{G}_a\times X_{m,n} &\rightarrow &X_{m,n} \\
(t,(x,y,u,v)) &\longmapsto &(x,y,u+tx^{m},v+ty^{n})
\end{eqnarray*}

The $\mathbb{G}_{a}$-action is generated by the derivation $\delta
:\,u\mapsto x^{m}\mapsto 0,\;v\mapsto y^{n}\mapsto 0$ and the so determined
principal $\mathbb{G}_{a}$-bundle $\pi :X_{m,n}\rightarrow \mathbb{A}_{\ast
}^{2}$ is trivialized over the open covering of $\mathbb{A}_{\ast }^{2}$ by
the open subsets 
\begin{equation*}
U_{x}=\{x\neq 0\}\text{ }\text{and }U_{y}=\{y\neq 0\}
\end{equation*}%
Writing $\theta :(\mathbb{G}_{m}\times \mathbb{G}_{a})\times
X_{m,n}\rightarrow X_{m,n}$ for the composite action, note that 
\begin{eqnarray*}
\theta ((\lambda ,t),(x,y,u,v)) &=&(\lambda x,\lambda y,\lambda
^{-n}u+t\lambda ^{-n}x^{m},\lambda ^{-m}v+t\lambda ^{-m}y^{n}) \\
&=&(\lambda x,\lambda y,\lambda ^{-n}u+t\lambda ^{-(m+n)}(\lambda
x)^{m},\lambda ^{-m}v+t\lambda ^{-(m+n)}(\lambda y)^{n}) \\
&=&\theta ((\lambda ,\lambda ^{-(m+n)}t),(x,y,u,v))
\end{eqnarray*}%
Set $d=m+n.$ From the last equation it is clear that $X_{m,n}$ is equipped
with a free action of the semidirect product $G_{d}=\mathbb{G}_{m}\ltimes
_{d}\mathbb{G}_{a}$ with multiplication 
\begin{equation*}
(\lambda ,t)\ (\lambda ^{\prime },t^{\prime })=(\lambda \lambda ^{\prime
},t+\lambda ^{d}t^{\prime }).
\end{equation*}

The open covering $\pi ^{-1}(U_{x})\cup $ $\pi ^{-1}(U_{x})$ that
trivializes the $\mathbb{G}_a$-action on $X_{m,n}$ is clearly $G_{d}$-stable
as well. We let

\begin{eqnarray*}
V_{x} &=&\pi ^{-1}(U_{x})\cong \mathrm{Spec}(\mathbb{C}[x^{\pm
1},y,u,v]/(x^{m}v-y^{n}u-1))\cong \mathrm{Spec}(\mathbb{C}[x^{\pm 1},y,u]) \\
V_{y} &=&\pi ^{-1}(U_{y})\cong \mathrm{Spec}(\mathbb{C}[x,y^{\pm
1},u,v]/(x^{m}v-y^{n}u-1))\cong \mathrm{Spec}(\mathbb{C}[x,y^{\pm 1},v]).
\end{eqnarray*}%
Setting 
\begin{equation*}
(u_{1},T_{1},L_{1})=(yx^{-1},x^{n}u,x)\text{ and }%
(u_{2},T_{2},L_{2})=(xy^{-1},y^{m}v,y)
\end{equation*}%
and observing that $u_{i}$ are $G_{d}$-invariant, $T_{i}$ are $\mathbb{G}%
_{m} $-invariant and translated by the $\mathbb{G}_{a}$-action (vice versa
for $L_{i}$), we obtain $G_{d}$-equivariant trivializations

\begin{eqnarray*}
V_{x} &\cong &\mathrm{Spec}(\mathbb{C}[x^{\pm 1},y,u])\cong \mathrm{Spec}(%
\mathbb{C}[u_{1}][T_{1},L_{1},L_{1}^{-1}])\cong \mathbb{A}^{1}\times G_{d} \\
V_{y}\ &\cong &\mathrm{Spec}(\mathbb{C}[x,y^{\pm 1},v\ ])\cong \mathrm{Spec}(%
\mathbb{C}[u_{2}][T_{2},L_{2},L_{2}^{-1}])\cong \mathbb{A}^{1}\times G_{d}.
\end{eqnarray*}

Gluing $V_{x}/G_{d}\cong \mathrm{Spec}(\mathbb{C}[u_{1}])$ and $%
V_{y}/G_{d}\cong \mathrm{Spec}(\mathbb{C}[u_{2}])$ over their intersection $%
\mathrm{Spec}(\mathbb{C}[u_{1},u_{1}^{-1}])$ via $u_{1}\mapsto u_{2}^{-1}$
we obtain a quotient 
\begin{equation*}
\pi _{m,n}:X_{m,n}\rightarrow C_{m,n}\cong \mathbb{P}^{1}
\end{equation*}%
which is a Zariski locally trivial principal $G_{d}$-bundle. Transition
isomorphisms are given by $L_{2}=y=u_{1}L_{1}$ and 
\begin{eqnarray*}
T_{2} &=&y^{m}v=(\frac{y}{x}%
)^{m}x^{m}v=u_{1}^{m}(1+y^{n}u)=u_{1}^{m}(1+u_{1}^{n}T_{1}) \\
&=&u_{1}^{m}+u_{1}^{d}T_{1}
\end{eqnarray*}%
proving

\begin{proposition}
The threefold $X_{m,n}$ admits the structure of a principal $G_{m+n}$-bundle
over $\mathbb{P}^{1}$ that is locally trivial for the Zariski topology. The
transition isomorphisms for this bundle are 
\begin{equation*}
(L,T)\longmapsto (uL,u^{m+n}T+u^{m}).
\end{equation*}
\end{proposition}

In contrast to the differential-topological context of \cite{tD} we cannot
conclude that the induced mapping $X_{m,n}/\mathbb{G}_m\rightarrow
C_{m,n}\cong \mathbb{P}^{1}$ is \emph{algebraically} isomorphic to the line
bundle $O_{\mathbb{P}^{1}}(-(m+n))$ (the translation part of the transition
isomorphism cannot be algebraically contracted away). Nevertheless, we can
conclude from the discussion above with $d=m+n,$ that

\begin{corollary}
\begin{enumerate}
\item The quotient $X_{m,n}/\mathbb{G}_m$ is isomorphic to the affine
surface $S_{d,m}$ which is the total space of the $\mathcal{O}_{\mathbb{P}%
^{1}}(-d)$-torsor $\rho :S_{d,m}\rightarrow \mathbb{P}^{1}$ with transition
isomorphism $T\mapsto u^{d}T+u^{m}$.

\item $\mathrm{Pic}(S_{d,m})\cong \mathrm{H}^{1}(S_{d,m},\mathcal{O}%
_{S_{d,m}}^{\ast })\cong \mathbb{Z}$ and the isomorphism class of the $%
\mathbb{G}_{m}$-bundle $X_{m,n}\rightarrow S_{d,m}$ is a generator of $%
\mathrm{Pic}(S_{d,m})$.
\end{enumerate}
\end{corollary}

\begin{proof}
The second assertion follows from the fact that every $\mathbb{A}^{1}$%
-bundle over $\mathbb{P}^{1}$ can be realized as the complement of a section 
$C$ in a suitable $\mathbb{P}^{1}$-bundle over $\mathbb{P}^{1}$ \cite[%
Theorem 2.3]{W} ($S_{d,m}$ is explicitly realized this way in the proof of
Theorem 1 below). The latter have divisor class group isomorphic to $\mathbb{%
Z}^{2},$ generated by the classes of a fiber and a section that can be
chosen to be precisely $C$. Thus the divisor class group of $S_{d,m}$ is
isomorphic to $\mathbb{Z}$, generated for instance by the class of a fiber
of $\rho $.
\end{proof}

The isomorphism class of $X_{m,n}$ as a variety is therefore determined by a
generator of \textrm{Pic}$(S_{d,m}),$ and moreover is independent of the
choice of generator. \ Indeed, choosing the other generator yields the $%
\mathbb{G}_m $-bundle $\tilde{X}_{m,n}$ over $S_{d,m}$ with transition
isomorphism $(L,T)\longmapsto (u^{-1}L,u^{d}T+u^{m})$. The coordinate change 
$\tilde{L}=L^{-1}$ on the fibers of the bundle yields an isomorphism $\tilde{%
X}_{m,n}\overset{\sim}{\rightarrow} X_{m,n}.$

The structure of $S_{d,m}$ as an $\mathcal{O}_{\mathbb{P}^{1}}(-d)$-torsor
(actually as an $\mathbb{A}^{1}$-bundle over $\mathbb{P}^{1})$ enables the
application of a famous theorem of Danilov-Gizatullin \cite{D-G} to
determine its isomorphy class as an affine surface (for a recent proof of
this theorem see \cite{fkz}). \ 

The transition isomorphism $\ T\longmapsto u^{d}T+u^{m}$ determines an $%
\mathbb{A}^{1}$-bundle over $\mathbb{P}^{1}$ which in turn can be completed
to the $\mathbb{P}^{1}$-bundle $\mathbb{P}(\mathcal{E)},$ where $\mathcal{E}$
is the rank 2 vector bundle over $\mathbb{P}^{1}$ with transition matrix 
\begin{equation*}
M = \left[ 
\begin{array}{cc}
u^{d} & u^{m} \\ 
0 & 1%
\end{array}
\right]
\end{equation*}
corresponding to a non trivial extension 
\begin{equation*}
0\rightarrow O_{\mathbb{P}^{1}}(-d)\rightarrow \mathcal{E}\rightarrow O_{%
\mathbb{P}^{1}}\rightarrow 0.
\end{equation*}%
Since the goal is to show that $X_{m,n}\cong X_{p,q}$ if $p+q=m+n,$ and
obviously $X_{m,n}\cong X_{n,m},$ we will assume that $m\geq n.$ The results
and terminology about ruled surfaces in what follows can be found in \cite[%
V.2]{H}.

\begin{lemma}
The total space of the $\mathbb{P}^{1}$-bundle $\mathbb{P}(\mathcal{\
E)\rightarrow }$ $\mathbb{P}^{1}$ is isomorphic to the Nagata-Hirzebruch
surface $\mathbb{F}_{2m-d}$.
\end{lemma}

\begin{proof}
It suffices to show that $\mathcal{E\otimes }\mathcal{O}_{\mathbb{P}^{1}}(m)$
is normalized, i.e. that%
\begin{equation*}
\mathrm{H}^{0}(\mathcal{E}\otimes \mathcal{O}_{\mathbb{P}^{1}}(m))\neq 0%
\text{ and }\mathrm{H}^{0}(\mathcal{E}\otimes \mathcal{O}_{\mathbb{P}%
^{1}}(m-1))=0.
\end{equation*}%
Indeed, if these conditions hold, then by counting degrees, $\mathcal{E}%
\otimes \mathcal{O}_{\mathbb{P}^{1}}(m)\cong \mathcal{O}_{\mathbb{P}%
^{1}}\oplus \mathcal{O}_{\mathbb{P}^{1}}(2m-d)$ and, as is well known, $%
\mathbb{P}(\mathcal{E)\cong }\mathbb{P}(\mathcal{E}\otimes \mathcal{O}_{%
\mathbb{P}^{1}}(m))\cong \mathbb{F}_{2m-d}$.

The vector bundle $\mathcal{E}\otimes \mathcal{O}_{\mathbb{P}^1}(j)$ has
transition matrix 
\begin{equation*}
M_{j}= \left[ 
\begin{array}{cc}
u^{d-j} & u^{m-j} \\ 
0 & u^{-j}%
\end{array}%
\right]
\end{equation*}%
and a nonzero global section for this bundle corresponds to a pair of
nonzero vectors 
\begin{eqnarray*}
g=\left[ 
\begin{array}{c}
g_{1} \\ 
g_{2}%
\end{array}
\right] \in \mathbb{C}[u]\times \mathbb{C}[u] & \text{ and } & h=\left[ 
\begin{array}{c}
h_{1} \\ 
h_{2}%
\end{array}
\right] \in \mathbb{C}[u^{-1}]\times \mathbb{C}[u^{-1}]
\end{eqnarray*}
with $M_{j}\cdot g=h$ For $j=m-1$ this would mean 
\begin{eqnarray*}
h_{1} &=&u^{d-m+1}g_{1}+ug_{2} \\
h_{2} &=&u(u^{d-m}g_{1}+g_{2})
\end{eqnarray*}%
which is impossible as $d>m$. \ For $j=m$ we have, for example, 
\begin{eqnarray*}
g=\left[ 
\begin{array}{c}
0 \\ 
1%
\end{array}%
\right] & \text{ and } & h=\left[%
\begin{array}{c}
1 \\ 
u^{-m}%
\end{array}%
\right].
\end{eqnarray*}
\end{proof}

\begin{theorem}
Let $d=p+q=m+n.$ Then $X_{m,n}\cong X_{p,q}$ as abstract varieties.
\end{theorem}

\begin{proof}
As noted above, since $X_{m,n}\cong X_{n,m},$ we may assume without loss of
generality that $m\geq n$, so that $2m\geq d.$ With the notation of the
corollary, $X_{m,n}\ $is the total space of a $\mathbb{G}_{m}$-bundle over $%
S_{d,m}$ which in turn is the total space of an $O_{\mathbb{P}^{1}}(-d)$%
-torsor over $\mathbb{P}^{1}$. Choosing homogeneous coordinates $[u:v]$ on
the fibers $F$ of the $\mathbb{P}^{1}$-bundle $\mathbb{P}(\mathcal{E)\cong 
\mathbb{F}}_{2m-d}\mathcal{\rightarrow }\mathbb{P}^{1}$, $S_{d,m}$ is
obtained as the complement of the section $v=0.$ The associated divisor is $%
C+mF$ where $C$ is the special section with self-intersection $d-2m$. The
self-intersection of $C+mF$ is calculated as 
\begin{equation*}
(C+mF)^{2}=C^{2}+2mC.F=-(2m-d)+2m=d
\end{equation*}%
and thus depends only on the sum $d=p+q=m+n.$ \ 

The theorem of Danilov-Gizatullin says exactly that $S_{d,m}$ and $S_{d,p}$
are isomorphic as affine surfaces. \ An isomorphism $S_{d,m}\overset{\phi }{%
\rightarrow }$ $S_{d,p}$ carries a generator of $\mathrm{Pic}(S_{d,p})$ to a
generator of $\mathrm{Pic}(S_{d,m})$ under the induced isomorphism$\ \ \ $%
\begin{equation*}
\mathrm{Pic}(S_{d,p})\overset{\phi ^{\ast }}{\rightarrow }\mathrm{Pic}%
(S_{d,m}).
\end{equation*}%
Since the class of $X_{m,n}$ (and $\tilde{X}_{m,n})$ generates $\mathrm{Pic}%
(S_{d,m})$ and $X_{m,n}$ $\cong \tilde{X}_{m,n}$ we obtain the isomorphism $%
X_{m,n}\cong X_{p,q}.$
\end{proof}

While it appears to be rather difficult to construct an explicit isomorphism
between even the first interesting examples, $X_{2,2}$ and $X_{3,1},$ one
can check that $X_{2,2}$ admits the structure of the principal $\mathbb{G}_a$%
-bundle over $\mathbb{A}^{2}_*$ corresponding to the \v{C}ech cocycle $%
x^{-3}y$.

\begin{example}
Set 
\begin{eqnarray*}
a(x,y,u,v) &=&x-\frac{1}{2}y \\
b(x,y,u,v) &=&\frac{(6x-y)\ }{8}v-\frac{(3y-2x)}{2}u \\
w(x,y,u,v) &=&\frac{5}{16}v^{2}x+\frac{5}{2}vxu-\frac{1}{32}v^{2}y-\frac{5}{4%
}vyu+u^{2}x-\frac{5}{2}u^{2}y
\end{eqnarray*}%
The morphism 
\begin{eqnarray*}
\tilde{\pi}\ &:&X_{2,2}\rightarrow \mathbb{A}^{2}_*\cong \mathrm{Spec}(%
\mathbb{C}[a,b])\setminus\! \{0\}. \\
& &(x,y,u,v)\mapsto (a(x,y,u,v),b(x,y,u,v))
\end{eqnarray*}%
is a $\mathbb{G}_a$-bundle corresponding to the \v{C}ech cocycle $a^{-3}b$.
In particular, one can check that the assignment 
\begin{equation*}
\ x\mapsto \frac{a^{3}}{3},\text{ }y\mapsto a^{3},\text{ }u\mapsto xb-\frac{1%
}{4},\text{ }v\mapsto 2yb-1
\end{equation*}%
extends to a locally nilpotent $\mathbb{C}$-derivation $\delta $ on $A_{2,2}=%
\mathbb{C}[x,y,u,v]/(x^{2}v-y^{2}u-1)$ with $a,b\in \ker (\delta )$ and 
\begin{equation*}
\delta (y+a+ab)=a^{3},\text{ }\delta (w)=b.
\end{equation*}%
Since $(a^{3},b)A_{2,2}=A_{2,2},$ we obtain a locally trivial $\mathbb{G}_a$%
-action on $X_{2,2}$ with quotient isomorphic to $\mathbb{A}^{2}_*$
corresponding to the \v{C}ech cocycle 
\begin{equation*}
\frac{y+a+ab}{a^{3}}-\frac{w}{b}=\frac{1}{a^{3}b}.
\end{equation*}
\end{example}

\section{A more general class of examples}

A refinement of the arguments in the previous section enables a description
of the isomorphy classes of total spaces of principal $\mathbb{G}_a$-bundles
defined by arbitrary homogeneous bivariate polynomials in terms of the $%
X_{m,n}$. \ Let $f,g\in \mathbb{C}\left[ x,y\right] $ be arbitrary
homogeneous polynomials such that $\mathrm{Spec}\left( \mathbb{C}\left[ x,y%
\right] /\left( f,g\right) \right)$ is supported at the origin of $\mathbb{A}%
^2$ and consider affine threefolds $X_{f,g}$ in $\mathbb{A}^{4}=\mathrm{Spec}%
\left( \mathbb{C} \left[ x,y,u,v\right] \right) $ defined by the equations $%
fv-gu-1=0$. Again, these threefolds come naturally equiped with the
structure $\pi _{f,g}:X_{f,g}\rightarrow \mathbb{A}^{2}_*$ of $\mathbb{G}%
_{a} $-bundles over the complement of the origin in $\mathbb{A}^2=\mathrm{%
Spec}\left( \mathbb{C}\left[ x,y\right] \right) $ by restricting the
projection $\mathrm{p}_{x,y}:\mathbb{A}^{4}\rightarrow \mathbb{A}^{2}$ to $%
X_{f,g}$. The latter becomes trivial on the covering of $\mathbb{A}^{2}_*$
by means of the principal affine open subsets 
\begin{equation*}
U_{f}=\mathrm{Spec}\left( \mathbb{C}\left[ x,y\right] _{f}\right) \text{ and 
}U_{g}=\mathrm{Spec}\left( \mathbb{C}\left[ x,y\right] _{g}\right) .
\end{equation*}

Suppose that $f$ and $g$ are homogeneous, say of degree $m$ and $n$
respectively. Similar to the case of the $X_{m,n}$ of the previous section,
the $\mathbb{G}_{m}$-action $\lambda \cdot \left( x,y\right) =\left( \lambda
x,\lambda y\right) $ on $\mathbb{A}^{2}$ lifts to a free $\mathbb{G}_{m}$
-action on $X_{f,g}$ defined by 
\begin{equation*}
\lambda \cdot \left( x,y,u,v\right) =\left( \lambda x,\lambda y,\lambda
^{-m}u,\lambda ^{-n}v\right)
\end{equation*}
which is compatible with the structural map $\pi _{f,g}:X_{f,g}\rightarrow 
\mathbb{A}^{2}_*$. It follows that $\pi _{f,g}$ descends to an $\mathbb{A}%
^{1}$-bundle $\rho _{f,g}:X_{f,g}/\mathbb{G}_{m}\rightarrow \mathbb{A}^2_* /%
\mathbb{G}_{m}\simeq \mathbb{P}^{1}$. The latter is an $\mathcal{O}_{\mathbb{%
P}^{1}}\left( -n-m\right) $-torsor which can be described more explicitly as
follows. Recall that for every pair of integers $m,n\geq 1$ , the quotient
of $\mathbb{A}^{2}_*\times \mathbb{A}^{2}_*$ by the free $\mathbb{G}_{m}^{2}$%
-action given by 
\begin{equation*}
\left( \lambda ,\mu \right) \cdot \left( x,y,u,v\right) =\left( \lambda
x,\lambda y,\mu \lambda ^{-n}u,\mu \lambda ^{-m}v\right)
\end{equation*}
is the surface scroll 
\begin{equation*}
\pi :\mathbb{F}\left( m,n\right) =\mathbb{P}\left( \mathcal{O}_{\mathbb{P}%
^{1}}\left( m\right) \oplus \mathcal{O}_{\mathbb{P}^{1}}\left( n\right)
\right) \rightarrow \mathbb{P}^{1},
\end{equation*}
with structural morphism $\pi $ induced by the $\mathbb{G}_{m}^{2}$%
-invariant morphism $\mathbb{A}^{2}_*\times \mathbb{A}^{2}_*\rightarrow 
\mathbb{P}^{1}$, $\left( x,y,u,v\right) \mapsto \left[ x:y\right]$. Clearly, 
$X_{f,g}$ is contained in the complement of the $\mathbb{G}_{m}^{2}$%
-invariant closed subset $D$ of $\mathbb{A}_{\ast }^{2}\times \mathbb{A}%
_{\ast }^{2}$ with equation $fv-gu=0$. Furthermore, $X_{f,g}$ is invariant
under the action of the first factor of $\mathbb{G}_{m}^{2}$ and is a slice
for the action of the second factor of $\mathbb{G}_{m}^{2}$ of $\mathbb{A}%
_{\ast }^{2}\times \mathbb{A}_{\ast }^{2}\setminus D$. Thus the quotient map
restricts to a $\mathbb{G}_{m}=\mathbb{G}_{m}\times \left\{ 1\right\} $%
-bundle 
\begin{equation*}
\rho :X_{f,g}\rightarrow X_{f,g}/\mathbb{G}_{m}\times \left\{ 1\right\}
\simeq \left( \mathbb{A}_{\ast }^{2}\times \mathbb{A}_{\ast }^{2}\setminus
D\right) /\mathbb{G}_{m}^{2}=\mathbb{F}\left( m,n\right) \setminus \Delta
_{f,g}.
\end{equation*}%
over the complement of the image $\Delta _{f,g}$ of $D$ in $\mathbb{F}(m-n)$.

The self-intersection of the curve $\Delta _{f,g},$ which is a section of $%
\pi$, can determined as follows. The images by the quotient map of the $%
\mathbb{G}_{m}^{2}$-invariant subsets $\left\{ v=0\right\} $ and $\left\{
u=0\right\} $ define two sections $C_{v}$ and $C_{u}$ of $\pi $ such that $%
C_{v}\cdot C_{u}=0$. Letting $L$ be a fiber of $\pi $, the divisors $%
mL+C_{v} $ and $nL+C_{u}$ are linearly equivalent since they differ by the
divisor of the rational function $h=fv/gu$ on $\mathbb{F}\left( m,n\right) $%
. Since $\Delta _{f,g}$ belongs to the complete linear system generated by
these divisors, we conclude that 
\begin{equation*}
\Delta _{f,g}^{2}=\left( mL+C_{v}\right) \cdot \left( nL+C_{u}\right) =m+n.
\end{equation*}

The above description implies that $X_{f,g}/\mathbb{G}_{m}$ is isomorphic to
the complement of a section with self-intersection $m+n$ in $\mathbb{F}%
\left( m,n\right) \simeq \mathbb{F}_{\left\vert m-n\right\vert }$, and so,
the $\mathbb{P}^{1}$-bundle $\mathbb{F}_{\left\vert m-n\right\vert
}\rightarrow \mathbb{P}^{1}$ restricts to an $\mathcal{O}_{\mathbb{P}%
^{1}}\left( -m-n\right) $-torsor 
\begin{equation*}
\rho _{f,g}:X_{f,g}/\mathbb{G}_{m}\rightarrow \mathbb{A}^2_*/\mathbb{G}%
_{m}\simeq \mathbb{P}^{1}.
\end{equation*}
Combined with the results of the previous section, this the leads to the
following

\begin{proposition}
With the notation above, $X_{f,g}$ is isomorphic to $X_{m,n}$.
\end{proposition}

\begin{proof}
Recall that $X_{m,n}$ has the structure of a $\mathbb{G}_{m}$-bundle over an
affine surface $S_{d,m}$ obtained as the complement of certain section $%
\Delta _{m,n}$ with self-intersection $m+n$ in $\mathbb{F}_{\left\vert
m-n\right\vert }$. By virtue of Gizatullin's classification \cite{D-G}, $%
\mathbb{F}_{\left\vert m-n\right\vert }\setminus \Delta _{f,g}\simeq \mathbb{%
F}_{\left\vert m-n\right\vert }\setminus \Delta _{m,n}=S_{d,m}$ as abstract
affine surfaces. It follows that $Z=X_{m,n}\times _{S_{d,m}}X_{f,g}$ is a $%
\mathbb{G}_{m}$-bundle over both $X_{m,n}$ and $X_{f,g}$ via the first and
the second projection respectively. Since the Picard groups of $X_{m,n}$ and 
$X_{f,g}$ are both trivial, these $\mathbb{G}_{m}$-bundles are actually
trivial, which yields an isomorphism $X_{m,n}\times \mathbb{A}_{\ast
}^{1}\simeq X_{f,g}\times \mathbb{A}_{\ast }^{1}$ as affine varieties. Since 
$X_{m,n}$ and $X_{f,g}$ are both irreducible and have no invertible
functions except nonzero constants, the latter descends to an isomorphism $%
X_{m,n}\simeq X_{f,g}$.
\end{proof}

\end{document}